\documentclass[12pt]{amsart}

\numberwithin{equation}{section}
\newcommand{\CC}{\mathbb{C}}

\newcommand{\cB}{\mathcal{B}}

\newcommand{\cL}{\mathcal{L}}

\newcommand{\cH}{\mathcal{H}}
\newcommand{\bL}{\mathbf{L}}

\DeclareMathOperator{\spa}{span}
\DeclareMathOperator{\col}{colspan}
\title[Realization of rational functions and applications]
{A new realization of rational functions, with applications to 
linear combination interpolation}
\usepackage{amsthm}
\usepackage{amsmath,mathrsfs}
\usepackage{color}
\usepackage{amssymb}
\usepackage{hyperref}
\usepackage{enumerate}
\usepackage{latexsym,array}
\usepackage{amsfonts}
\usepackage{shadow}
\newtheorem{Pa}{Paper}[section]
\newtheorem{Tm}[Pa]{{\bf Theorem}}

\newtheorem{Cy}[Pa]{{\bf Corollary}}

\newtheorem{Rk}[Pa]{{\bf Remark}}
\newtheorem{Pn}[Pa]{{\bf Proposition}}
\newtheorem{Pb}[Pa]{{\bf Problem}}

\author[D. Alpay]{Daniel Alpay}
\address{(DA) Department of Mathematics\\
Ben-Gurion University of the Negev\\
Beer-Sheva 84105 Israel} \email{dany@math.bgu.ac.il}

\author[P. Jorgensen]{Palle Jorgensen}
\address{(PJ)  Department of Mathematics\\
 University of Iowa.
Iowa City, IA 52242 USA}
\email{palle-jorgensen@uiowa.edu}

\author[I. Lewkowicz]{Izchak Lewkowicz}
\address{(IL) Department of Electrical \& Computer Engineering\\
Ben-Gurion University of the Negev\\
Beer-Sheva 84105 Israel} \email{izchak@ee.bgu.ac.il}

\author[D. Volok]{Dan Volok}
\address{(DV) Department of Mathematics\\ Kansas State University\\
Manhattan KS 66506 USA
}
\email{danvolok@math.ksu.edu}


\thanks{The authors thank the Binational
Science Foundation Grant number 2010117. D. Alpay thanks the Earl
Katz family for endowing the chair which supported his research.}
\begin{document}
\parindent 0cm

\parindent 0cm
\begin{abstract}
We introduce the following {\sl linear combination interpolation} problem
(LCI), which in case of simple nodes reads as follows: Given $N$  distinct numbers $w_1,\ldots w_N$ and $N+1$ complex numbers $a_1,\ldots, a_N$ and $c$, find all 
functions $f(z)$ analytic in an open set (depending on $f$) containing the points $w_1,\ldots,w_N$ such that
\[
\sum_{u=1}^Na_uf(w_u)=c.
\]
To this end we prove a representation theorem for such functions $f$ in terms of an associated
polynomial $p(z)$. We give applications of this representation theorem to realization of rational functions and representations of positive definite kernels.
\end{abstract}
\maketitle
\tableofcontents
\section{Introduction}
\setcounter{equation}{0}

Any function $f$ analytic in a neighborhood of the origin can be uniquely written as 
\begin{equation}
\label{decomp1}
f(z)=\sum_{n=0}^{N-1}z^nf_n(z^N),
\end{equation}
where $f_0,\ldots, f_{N-1}$ are analytic at the origin. Furthermore, the maps 
\begin{equation}
\label{iowa-city-180215}
T_nf=f_n,\quad n=0,\ldots,N-1
\end{equation}
satisfy, under appropriate hypothesis, the Cuntz relations (see also \eqref{explicit}). See for instance \cite{MR3050315}, where applications to wavelets are given.
In the present paper, we extend these methods, and we derive new and explicit formulas for solutions to a class of multi-point interpolation problems; not amenable to 
tools from earlier investigations. Following common use, by “Cuntz relations” we refer here to a symbolic representation of a finite set (say  $N$) of isometries having 
orthogonal ranges which add up to the identity (operator). When $N$ is fixed the notation  $O_N$ is often used. By a representation of $O_N$ in a fixed Hilbert space, 
we mean a realization of the $N$-Cuntz relations in a Hilbert space. Here we will be applying this to specific Hilbert spaces of analytic functions which are dictated by our 
multi-point interpolation setting. 
In general it is known that the problem of 
finding representations of $O_N$ is subtle. (The literature on representations of $O_N$ is vast.) 
For example, no complete classification of these representations is known, but nonetheless, specific representations can be found, and they are known to play a key role in 
several areas of mathematics and its applications; e.g., to multi-variable operator theory, and in applications to the study of multi-frequency bands; see the cited references 
below. Realizations as in \eqref{decomp1} then results from representations of 
$O_N$; the particulars of these representations are then encoded in the operators from \eqref{iowa-city-180215}. Readers not familiar with $O_N$ and its representations are 
referred to \cite{BrJo02a, MR1284951,MR2254502,MR2097020}, and to Remark \ref{rk35} below.\\

The outline of the paper is as follows. In Section 2, we replace $z^N$ by an arbitrary polynomial $p(z)$, and prove a counterpart of the decomposition \eqref{decomp1}, see Theorem \ref{fbdthm}. 
The rest of the paper is organized as follows: In sections 3 -4, we discuss uniqueness of solutions, and (motivated by applications from systems theory) we extend our result in three ways, first  to that of Banach space valued functions (Theorem \ref{vvf}) and then we specialize to the case rational functions (Theorem \ref{realrat}). Thirdly we study multipoint interpolation when derivatives are specified. In section 5, we give an application to positive definite kernels.\\

More precisely, a first application of Theorem \ref{fbdthm} is in giving a new realization formula for rational functions. To explain the result, recall that a matrix-valued rational function $W$ analytic at the origin can always be written in the form
\[
W(z)=D+zC(I-zA)^{-1}B,
\]
where $D=W(0)$ and $A,B,C$ are matrices of appropriate sizes.
Such an expression is called a state space realization, and plays an important role in linear system theory and related topics; see
for instance \cite{bgr} and \cite{MR2663312} for more information.
We here prove that a rational function analytic at the points $w_1,\ldots, w_n$ can always be written in the form
\[
W(z)=Z(z)\gamma(I-p(z)\alpha)^{-1}\beta,
\]
where $p(z)$ is a polynomial vanishing at the points $w_1,\ldots, w_n$ and of degree $N\ge n$, and
\begin{equation}
Z(z)=\begin{pmatrix} 1&z&\cdots&z^{N-1}\end{pmatrix},
\label{Zz}
\end{equation}
and $\alpha, \beta,\gamma$ are matrices of appropriate sizes.
Finally we give an application to decompositions of positive definite kernels and the Cuntz relations.\smallskip

The multipoint interpolation problem,
which in the case where $p$ has simple zeros $w_1,\ldots,w_N$ consists in finding all 
functions $f(z)$ analytic in a simply connected set (depending on $f$) containing the points $w_1,\ldots,w_N$ and such that
\begin{equation}
\label{cons1}
\sum_{u=1}^Na_uf(w_u)=c.
\end{equation}

This can be equivalently written as
\[
(a_1,~\ldots~,~a_N)\left(\begin{smallmatrix}f(w_1)\\ \vdots \\ f(w_N)\end{smallmatrix}\right)=c.
\]
Namely the points $f(w_1),~\ldots~,~f(w_N)$ lie on a hyper-plane, so roughly speaking, the
points $w_1,~\ldots~,~w_N$ lie on some manifold.\smallskip

This type of problem seems to have been virtually neglected
in the litterature. In \cite{MR99h:47021} the case of two points was considered in the setting of the Hardy space of the open unit disk. The method there consisted in finding an involutive self-map of the open unit disk mapping one of the points to the second one, and thus reducing the
given two-point interpolation problem to a one-point interpolation problem with an added symmetry. This method cannot be extended to more than two points, but in special cases. In \cite{ajlm1} we considered the interpolation condition \eqref{cons1} in the Hardy space. Connections with the Cuntz relations played a key role in the arguments.\smallskip

\section{A decomposition of analytic functions}
\setcounter{equation}{0}
We set
$$p(z)=\prod_{j=1}^n(z-w_j)^{\mu_j},\quad \sum_{j=1}^n\mu_j=N,$$
and recall that $Z(z)$ is given by \eqref{Zz}.

\begin{Tm} \label{fbdthm}
Let $\Omega$ be a (possibly disconnected) neighborhood of $\{w_1,\dots, w_n\}$. Then there exists a neighborhood $\Omega_0$ of the origin, such that $ p^{-1}(\Omega_0)\subset\Omega$
 and every function
$f(z)$, analytic in $\Omega$, can be represented in the form
\begin{equation}\label{fbd}
 f(z)=Z(z)F(p(z)), \quad z\in p^{-1}(\Omega_0),
\end{equation}
where 
$F(z)$ is a $\CC^{N}$-valued function, analytic in $\Omega_0$. 
\end{Tm}

\begin{proof}
Choose $n$  simple closed counterclockwise oriented contours $\gamma_1,\dots,\gamma_n$ with the following properties:
\begin{enumerate}
 \item The function $f(z)$ is analytic on each contour $\gamma_j$ and in the simply connected domain $D_j$ encircled by $\gamma_j$.
\item For $j=1,\dots, n$ the domain $D_j$ contains the point $w_j$.
\item The domains $D_1,\dots,D_n$ are pairwise disjoint.
\end{enumerate}
 
Denote
$$D:=\bigcup_{j=1}^n D_j,\quad\rho:=\min\left\{|p(s)|\,:\,s\in\bigcup_{j=1}^n\gamma_j\right\}.$$
Since all the zeros of $p(z)$ are contained in $D$, $\rho>0$ and, by the maximum modulus principle, $p^{-1}(\Omega_0)\subset D,$ 
where $\Omega_0$ is  the open disk of radius $\rho$ centered at the origin. 
Furthermore,
for $z\in p^{-1}(\Omega_0)$ it holds that
\begin{multline*}
f(z)=\dfrac{1}{2\pi i }\sum_{j=1}^n\int_{\gamma_j}\dfrac {f(s)}{p(s)-p(z)}\dfrac{p(s)-p(z)}{s-z}ds
\\=\dfrac{Z(z)}{2\pi i }\sum_{j=1}^n\int_{\gamma_j}\dfrac {Q(s)f(s)}{p(s)-p(z)}ds=Z(z)F(p(z)),\end{multline*}
where $Q(s)$ is a $\CC^{N}$-valued polynomial, such that
$$\dfrac{p(s)-p(z)}{s-z}=Z(z)Q(s),$$
and
\begin{equation}
\label{eqF}
F(z)=\dfrac{1}{2\pi i }\sum_{j=1}^n\int_{\gamma_j}\dfrac {Q(s)f(s)}{p(s)-z}ds,\quad z\in \Omega_0,
\end{equation}
is a $\CC^{N}$-valued function, analytic in $\Omega_0$.
\end{proof}

\begin{Cy}
\label{cy123}
Assume that $f$ is a polynomial (resp. rational). Then $F$ given by \eqref{eqF} is also a polynomial (resp. rational).
\end{Cy}

\begin{proof}
We first consider the case of a polynomial. For $z$ near the origin we have
\[
F(z)=\sum_{u=0}^\infty z^uF_u\quad\text{with}\quad
F_u=\sum_{j=1}^n\frac{1}{2\pi i}\int_{\gamma_{j}}\frac{Q(s)f(s)}{p(s)^{u+1}}du.
\]
Note that $\frac{1}{2\pi i}\int_{\gamma_{j}}\frac{Q(s)f(s)}{p(s)^{u+1}}du$ is the residue of the rational function
$\frac{Q(s)f(s)}{p(s)^{u+1}}$ at the point $w_j$.
For $u$ large enough the difference of the degrees of the denominator and the numerator of this rational function is at least two, and so the 
sum of its residues is equal to $0$ (the so-called exactity relation; see \cite[p. 173]{freitag2} and \cite[Exercise 7.3.6, p. 326]{CAPB}). Thus $F_u=0$ for $u$ 
large enough and we conclude by analytic continuation that $F$ is a polynomial.\smallskip

In the case of a rational function consider  the partial fraction representation, which is the sum of
a polynomial (which we just have treated) and of terms of the form $\frac{1}{(s-a)^M}$, where $a$ is not a zero of $p$. We thus need to show that, for such $a$, a sum of the form
\begin{equation}
\label{eqFF}
G(z)=\dfrac{1}{2\pi i }\sum_{j=1}^n\int_{\gamma_j}\dfrac{Q(s)}{(s-a)^M(p(s)-z)}ds,
\end{equation}
is rational. Chose the contours $\gamma_1,\ldots, \gamma_n$ such that no zeroes of the equation $p(s)=p(a)$ lie inside or on them.
Using the polynomial case, the result follows from 
writing
\[
\begin{split}
G(z)&=\dfrac{1}{2\pi i }\sum_{j=1}^n\int_{\gamma_j}\dfrac{Q(s)}{(s-a)^M(p(s)-z)}ds\\
&=\dfrac{1}{2\pi i }\sum_{j=1}^n\int_{\gamma_j}\dfrac{Q(s)\left(\frac{p(s)-p(a)}{s-a}\right)^M}{(p(s)-p(a))^M(p(s)-z)}ds\\
&=\dfrac{1}{2\pi i }\sum_{j=1}^n\int_{\gamma_j}Q(s)\left(\frac{p(s)-p(a)}{s-a}\right)^M
\left(\frac{c(z)}{p(s)-z}+\right.\\
&\hspace{6cm}+\left. \sum_{u=1}^M\frac{c_u(z)}{(p(s)-p(a))^u}
\right)ds,
\end{split}
\]
for some  complex numbers $c(z),c_1(z),\ldots, c_M(z)$ corresponding to the partial fraction expansion of the function $\frac{1}{(\lambda-z)(\lambda-p(a))^M}$:
\[
\frac{1}{(\lambda-z)(\lambda-p(a))^M}=\frac{c(z)}{\lambda-z}+\sum_{u=1}^M\frac{c_u(z)}{(\lambda-p(a))^u}.
\]
These are readily seen to be rational functions of $z$, and hence the function $G$ above is rational.
\end{proof}

\section{A new realization of rational functions}
Denote by $V$ the generalized $N\times N$ Vandermonde matrix
$$V=\begin{pmatrix} V_1 \\ \vdots \\ V_n\end{pmatrix},\text{ where }\,
V_j=\begin{pmatrix} Z(w_j) \\ Z^\prime(w_j)\\ \vdots \\ Z^{(\mu_j-1)}(w_j)\end{pmatrix},$$
and by $C_w$ the linear operator
$$f(z)\mapsto\begin{pmatrix}C_{w_1}f\\ \vdots \\ C_{w_n}f\end{pmatrix},\text{ where }\,
C_{w_j}f:=\begin{pmatrix} f(w_j)\\ f^\prime(w_j)\\ \vdots \\ f^{(\mu_j-1)}(w_j)\end{pmatrix}.$$

By rearanging the rows the matrix $V$ is readily seen to be invertible.

\begin{Pn}
\label{pn31}
Let $f(z)$ be a function, analytic in a neighborhood of $\{w_1,\dots,w_n\}$ and let $F(z)$ be a $\CC^{N}$ function,
 analytic in a neighborhood of the origin, which provides the decomposition \eqref{fbd} for the function $f(z)$. Then the Taylor expansion 
of $F(z)$  is given by 
\begin{equation}
\label{u1}
F(z)=\sum_{k=0}^\infty z^kV^{-1}C_w(R_0^{(p)})^kf,
\end{equation}
where $R_0^{(p)}$ denotes the linear operator 
$$f(z)\mapsto\dfrac{f(z)-Z(z)V^{-1}C_wf}{p(z)}.$$
\end{Pn}

\begin{Rk}
{\rm 
{\it A priori} the convergence in \eqref{u1} is pointwise, and uniform on compact subsets of the origin where $f$ is defined. When the underlying spaces are finite dimensional, or when some extra topological structure is given, one can rewrite \eqref{u1} as
\[
F(z)=V^{-1}C_w(I-zR_0^{(p)})^{-1}f.
\]
}
\end{Rk}
\begin{proof}[Proof of Proposition \ref{pn31}]
Since 
$$p(w_j)=p^\prime(w_j)=\dots=p^{(\mu_j-1)}(w_j)=0,\quad j=1,\dots,n,$$
differentiate both sides of \eqref{fbd} at $w_j$ to obtain
$$C_{w_j}f=V_jF(0),\quad j=1,\dots,n.$$
Hence, in vector notation,
$$C_wf=VF(0),$$
\begin{equation}F(0)=V^{-1}C_wf,\label{u2}\end{equation}
and
\begin{equation}\label{u3}
 (R_0^{(p)}f)(z)=Z(z)(R_0F)(p(z)),
\end{equation}
where
$$(R_wF)(z)=\dfrac{F(z)-F(w)}{z-w}$$
is the classical backward-shift operator.
In particular, the function $R_0F$ provides a decomposition \eqref{fbd} for the function $R_0^{(p)}f$. By induction, one may conclude that
$$((R_0^{(p)})^kf)(z)=Z(z)R_0^kF(p(z)),\quad k=0,1,2,\dots$$
hence, in view of \eqref{u2},
$$(R_0^kF)(0)=V^{-1}C_w(R_0^{(p)})^kf,\quad k=0,1,2,\dots$$
\end{proof}
\begin{Cy}\label{refmist}
Every function $f(z)$, analytic in a neighborhood of  $\{w_1,\dots,w_n\}$, admits a unique decomposition \eqref{fbd}, in which (as follows form Corollary \ref{cy123})
$F$ is a polynomial (resp. rational) when $f$ is a polynomial (resp. rational).
\end{Cy}

Theorem \ref{fbdthm} has an analogue in the setting of analytic functions with values in a Banach space $\cB$.
In what follows, $\cB^s$ denotes the product  space 
$$\cB^s:=\CC^s\otimes\cB,$$
and the tensor product of a matrix $(a_{i,j})\in\CC^{r\times s}$ and a linear operator $A\in\cL(B)$ is understood as the operator matrix
$$(a_{i,j})\otimes A:=(a_{i,j}A)\in\cL(\cB^s,\cB^r).$$ 
\begin{Tm}\label{vvf}
Let $\cB$ be a Banach space and let $f(z)$ be a $\cB$-valued function, analytic in a neighborhood $\Omega$ of $\{w_1,\dots,w_n\}$. Then there exist a neighborhood $\Omega_0$
of the origin,  and a $\cB^N$-valued function $F(z)$, analytic in $\Omega_0$, 
such that
\begin{equation}\label{fbd2}f(z)=(Z(z)\otimes I_\cB)F(p(z)),\quad z\in p^{-1}(\Omega_0)\subset\Omega.\end{equation}
Furthermore, the Taylor expansion 
of $F(z)$  is given by 
\begin{equation}
F(z)=\sum_{k=0}^\infty z^kV^{-1}C_w(R_0^{(p)})^kf,\label{u101}
\end{equation}
where $R_0^{(p)}$ denotes the linear operator 
$$f(z)\mapsto\dfrac{f(z)-((Z(z)V^{-1})\otimes I_\cB)C_wf}{p(z)}.$$
\end{Tm}

\begin{proof}
Let $\varphi\in\cB^*$. Then, according to Theorem \ref{fbdthm} and Proposition \ref{pn31},
 the function $\varphi\circ f$ admits a unique decomposition \eqref{u1} provided by the $\CC^N$-valued function 
$$F^{\varphi}(z)=\sum_{k=0}^\infty z^kV^{-1}C_w(R_0^{(p)})^k(\varphi\circ f).$$

Then
$$F^{\varphi}(z)=\sum_{k=0}^\infty z^k(I_N\otimes\varphi)F_k,$$
where $$F_k:=(V^{-1}\otimes I_\cB)C_w(R_0^{(p)})^kf.$$
Since $F^{\varphi}(z)$ is analytic in an open disk  $$\Omega_0=\{z\,:\,|z|<\rho\},$$
where $\rho$ is independent of $\varphi$,  the uniform boundedness principle implies that the $\cB^N$-valued function
$$F(z):=\sum_{k=0}^\infty z^kF_k$$ is also analytic in $\Omega_0,$ and \eqref{fbd2} follows from
$$F^\varphi=(I_N\otimes\varphi)\circ F,\quad \varphi\in\cB^*.$$

\end{proof}

The preceding analysis leads to a new kind of realization for rational functions. 
\begin{Tm}
\label{realrat}
Every rational $\CC^{r\times s}$-valued function $f(z)$, which has no poles in $\{w_1,\dots,w_n\},$ can be written as
\begin{equation}
\label{newreal}
f(z)=(Z(z)\otimes I_r)C(I-p(z)A)^{-1}B,
\end{equation}
where $A,B,C$ are constant matrices of appropriate sizes.
\end{Tm}
\begin{proof}
Write
$$\dfrac{1}{p(z)}=\sum_{j=1}^n\sum_{k=1}^{\mu_j}\dfrac{c_{j,k}}{(z-w_j)^k},$$
where $c_{j,k}\in\CC$ are constants.
Then the operator $R_0^{(p)}$ defined in \eqref{u3} can be written as
$$R_0^{(p)}=\sum_{j=1}^n\sum_{k=1}^{\mu_j}c_{j,k}R_{w_j}^k.$$
Since $f(z)$ is a rational function, the space
\[
\begin{split}
\cL(f)&:=\col\{A^kf\,:\,k=0,1,2,\dots\}\\
&\subset\col\{R_{w_j}^kf\,:\,j=1,\dots,n;\ k=0,1,2,\dots\}
\end{split}
\]
is finite-dimensional. Choose a basis of this finite-dimensional space and let $A,B,C$ be matrices representing the operators
\begin{gather*}\cL(f)\ni fu\mapsto R_0^{(p)}fu\in\cL(f),\quad u\in\CC^s,\\
 \CC^s\ni u\mapsto fu\in\cL(f),\\
\cL(f)\ni fu\mapsto (V^{-1}\otimes I_r)C_wfu\in\CC^{rN}, \quad u\in\CC^s,\end{gather*} respectively. 
Then
$$f(z)=(Z(z)\otimes I_r)\sum_{k=0}^\infty p(z)^kCA^kB=(Z(z)\otimes I_r)C(I-p(z)A)^{-1}B.$$

\end{proof}

We will call the realization \eqref{newreal} minimal if the size of the matrix $A$ is minimal (for more on this notion, and equivalent characterizations, see for instance
\cite{bgk1}). As a consequence of the uniqueness of the decomposition we also have:

\begin{Cy}
When minimal, the realization \eqref{newreal} is unique up to a similarity matrix. Then, $F$ is a polynomial if and only if $A$ is nilpotent.
\end{Cy}

\begin{proof}
It suffices to notice that the uniqueness of the decomposition \eqref{fbd}  reduces the problem to the uniqueness of the minimality of the function 
$C(I-\lambda A)^{-1}B$ with $\lambda\in\mathbb C$.
\end{proof}
\begin{Rk}
{\rm
The uniqueness allows us to give an interpretation on terms of generalized Cuntz relations. More precisely, define linear operators on analytic functions by
$S_1,\ldots, S_N,T_1,\ldots, T_N$ by:
\begin{equation}
\label{explicit}
(S_jg)(z)=z^{j-1}g(p(z))\quad \text{and}\quad T_jF=F_j,\quad j=1,\ldots,N 
\end{equation}
where $F$ is a $\mathbb C^N$-valued analytic function (see also \eqref{iowa-city-180215} for the defintion of $T_1,\ldots,T_N$).Then the given decomposition \eqref{fbd} reads
\[
T_iS_j=\delta_{ij}\quad\text{and}\quad \sum_{n=1}^{N} S_jT_j=I.
\]
}
\label{rk35}
\end{Rk}

\begin{Rk} {\rm We note that in the case $\mu_1=\dots=\mu_n=1$ the operator $R_0^{(p)}$  can be written as 
\begin{equation}
\label{aff}
R_0^{(p)}f(z)=\frac{f(z)}{p(z)}-\sum_{u=1}^N\dfrac{f(w_u)}{p^\prime(w_u)(z-w_u)}
\end{equation}
is reminiscent of a formula for a resolvent operator given in the setting of function theory on compact real
Riemann surfaces. See \cite[(4.1), p. 307]{av3}. This point is emphasized in the following proposition.}
\end{Rk}

\begin{Pn}
Let $\alpha$ and $\beta$ be such that the roots $w_1(\alpha),\ldots, w_N(\alpha)$ and $w_1(\beta),\ldots, w_N(\beta)$
of the equations $p(z)=\alpha$ and $p(z)=\beta$ are all distinct ($w_u(\alpha)\not=w_v(\beta)$ for $u,v=1,\ldots,N$). 
Then the resolvent equation
\[
R^{(p)}_\alpha-R^{(p)}_\beta=(\alpha-\beta)R^{(p)}_\alpha R^{(p)}_\beta
\]
holds.
\label{vinnikov}
\end{Pn}

\begin{proof} Indeed, on the one hand,
\[
\begin{split}
\left((R^{(p)}_\alpha-R^{(p)}_\beta)(f)\right)(z)&=\\
&\hspace{-5cm}=\frac{f(z)}{p(z)-\alpha}-\sum_{u=1}^N\dfrac{f(w_u(\alpha))}{p^\prime(w_u(\alpha))(z-w_u(\alpha))}-
\frac{f(z)}{p(z)-\alpha}+\sum_{u=1}^N\dfrac{f(w_u(\alpha))}{p^\prime(w_u(\alpha))(z-w_u(\alpha))}\\
&\hspace{-5cm}
=(\alpha-\beta)\frac{f(z)}{(p(z)-\alpha)(p(z)-\beta)}-\sum_{u=1}^N\dfrac{f(w_u(\alpha))}{p^\prime(w_u(\alpha))(z-w_u(\alpha))}+
\sum_{v=1}^N\dfrac{f(w_v(\beta))}{p^\prime(w_v(\beta))(z-w_v(\beta))}.
\end{split}
\]
On the other hand,
\[
\begin{split}
\left((R^{(p)}_\alpha\left(R^{(p)}_\beta(f)\right)\right)(z)&=\frac{\left(R^{(p)}_\beta(f)\right)(z)}{p(z)-\alpha}-
\sum_{u=1}^N\frac{\left(R^{(p)}_\beta(f)\right)(w_u(\alpha))}{p^\prime(w_u(\alpha))(z-w_u(\alpha))}\\
&=\frac{f(z)}{(p(z)-\alpha)(p(z)-\beta)}-
\sum_{v=1}^N\frac{f(w_v(\beta))}{p^\prime(w_v(\beta))(z-w_v(\beta))(p(z)-\alpha)}-\\
&\hspace{5mm}-
\sum_{u=1}^N\frac{f(w_u(\alpha))}{(
\underbrace{p(w_u(\alpha)}_{\text{$=\alpha$}}-\beta))p^\prime(w_u(\alpha))(z-w_u(\alpha))}+\\
&\hspace{5mm}+\sum_{u=1}^N
\frac{\left(\sum_{v=1}^N\frac{f(w_v(\beta))}{p^\prime(w_v(\beta))(w_u(\alpha)-w_v(\beta))}\right)}{p^\prime(w_u(\alpha))(z-w_u(\alpha))}.
\end{split}
\]
Proving the resolvent identity amounts to showing that:
\begin{equation}
\label{tyui}
\begin{split}
\sum_{v=1}^N\dfrac{f(w_v(\beta))}{p^\prime(w_v(\beta))(z-w_v(\beta))}=\\
&\hspace{-5cm}=(\alpha-\beta)\left\{
-
\sum_{v=1}^N\frac{f(w_v(\beta))}{p^\prime(w_v(\beta))(z-w_v(\beta))(p(z)-\alpha)}+
\sum_{u=1}^N
\frac{\left(\sum_{v=1}^N\frac{f(w_v(\beta))}{p^\prime(w_v(\beta))(w_u(\alpha)-w_v(\beta))}\right)}{p^\prime(w_u(\alpha))(
z-w_u(\alpha))}
\right\}\\
&\hspace{-5cm}=(\alpha-\beta)\left\{
-
\sum_{v=1}^N\frac{f(w_v(\beta))}{p^\prime(w_v(\beta))(z-w_v(\beta))(p(z)-\alpha)}+\right.\\
&\hspace{-4.5cm}\left.
+\sum_{v=1}^N
\frac{f(w_v(\beta))}{p^\prime(w_v(\beta))}\left\{\sum_{u=1}^N\frac{1}{
p^\prime(w_u(\alpha))(z-w_u(\alpha))(z-w_v(\beta))}+\right.\right.\\
&\left.\left. +\frac{1}{p^\prime(w_u(\alpha))(w_u(\alpha)-w_v(\beta))(z-w_v(\beta))}\right\}\right\}.
\end{split}
\end{equation}
Taking into account the equality
\[
\frac{1}{p(z)-\alpha}=\sum_{u=1}^N\frac{1}{p^\prime(w_u(\alpha))(z-w_u(\alpha))}
\]
we have
\[
\sum_{u=1}^N\frac{1}{p^\prime(w_u(\alpha))(w_u(\alpha)-w_v(\beta))}=\frac{1}{\alpha-p(w_v(\beta))}=\frac{1}{\alpha-\beta}.
\]
\eqref{tyui} follows.
\end{proof}

\begin{Rk}{\rm
Proposition \ref{vinnikov} can be proved in an easier way using the classical resolvent identity by remarking that
$(R_\alpha^{(p)}f)(z)=Z(z)(R_\alpha F)(p(z))$, where $f(z)=Z(z)F(p(z))$. The proof proposed here is more conducive to explicit links
with the Riemann surface case.
}
\end{Rk}
We finally note that, at least in spirit, we used the theory of linear system in this section. See for instance 
\cite{MR0255260,MR2663312,MR2002b:47144} 
for more on this theory. In the sequel we resort to the theory of reproducing kernel Hilbert spaces.
 The reader may find the following references helpful: \cite{saitoh,meschkowski,MR2002b:47144,MBM+14, NZ13, FM13, Wu13, Mou12}

\section{Linear combination interpolation}

In \cite{ajlm1} we introduced a general problem of linear combination interpolation, and solved it in the setting of the Hardy space.
Here, the preceding analysis enables us to solve a linear combination interpolation problem in the setting of functions analytic in the neighborhoods of given preassigned points.
\begin{Pb}
\label{newpb}
Given complex numbers $a_{j,k}, j=1,\dots,n;k=0,\dots,\mu_j-1;$ and $c$, describe the set of all functions $f$ analytic in a possibly disconnected neighborhood
of the points $w_1,\ldots, w_n$ and such that
\begin{equation}
\label{newinter}
\sum_{j=1}^n\sum_{k=0}^{\mu_j-1}a_{j,k}f^{(k)}(w_j)=c.
\end{equation}
\end{Pb}
The idea is to use the decomposition \eqref{u1} and to reduce the interpolation condition \eqref{newinter} to a {\sl unique} interpolation condition 
for a {\sl vector-valued} analytic function. Let 
$$v=\begin{pmatrix}a_{1,0}& a_{1,1}&\cdots&a_{n,\mu_n-1}\end{pmatrix}.$$
Then \eqref{newinter} can be re-written as
$$vC_wf=c.$$
In view of Propositions \ref{pn31}, this last condition is equivalent
to
$$vVF(0)=c,$$
which is a basic interpolation problem whose solution is given by
\[
F(z)=\frac{V^*v^*}{vVV^*v^*}c+\left(I_N+(z-1)\frac{V^*v^*vV^*}{{vVV^*v^*}}G(z)\right),
\]
where $G(z)$ is an arbitrary $\mathbb C^N$-valued function analytic in a neighborhood of the origin. Thus the solutions $f$ are given by
\[
\begin{split}
f(z)&=Z(z)F(p(z))\\
&=Z(z)\left(
\frac{V^*v^*}{vVV^*v^*}c+\left(I_N+(p(z)-1)\frac{V^*v^*vV^*}{{vVV^*v^*}}G(p(z))\right)\right).
\end{split}
\]
Furthermore, we obtain all the rational solutions of the interpolation when $G(z)$ is chosen rational.

\section{Representation in reproducing kernel Hilbert spaces}
Here we focus on the case when $f(z)$ belongs to a reproducing kernel Hilbert space $\cH(K)$ of analytic functions.
\begin{Pn} Let $K(z,w)$ be a positive definite function analytic in $z$ an in $\overline{w}$ in an open set $\Omega$ which contains $w_1,\dots,w_n.$ There exists a neighborhood $\Omega_0$ of the origin and a positive $\CC^{N\times N}$-valued kernel $L(z,w)$, analytic in $\Omega_0$, such that
$$K(z,w)=Z(z)L(p(z),p(w))Z(w)^*,\quad z,w\in p^{-1}(\Omega_0)\subset\Omega.$$
\end{Pn}

\begin{proof} Write $K(z,w)=C(z)C(w)^*$, where $C(z):\cH(K)\longrightarrow\CC$ is the point evaluation functional:
$$C(z)f=f(z),\quad z\in\Omega.$$
Then $C(z)$ is a $\bL(\cH(K),\CC)$-valued function, analytic in $\Omega$ and, by Theorem \ref{vvf},
 there exists a neighborhood $\Omega_0$ of the origin, such that
$p^{-1}(\Omega_0)\subset\Omega$  and
 $$C(z)=Z(z)E(p(z)),\quad z\in\Omega_0,$$ where
$E(z)$ is a $\bL(\cH(K),\CC^{N\times 1})$-valued function, analytic in $\Omega_0.$
Now set
$$L(z,w):=E(z)E(w)^*$$ to compete the proof.
\end{proof}

\begin{Pn}\label{compmult} Let $F\in\cH(L)$. Then the function
$Z(z)F(p(z)),$ which is analytic \textit{a priori} in $p^{-1}(\Omega_0)$, admits analytic continuation into $\Omega$ and is an 
element of the reproducing kernel Hilbert space $\cH(K).$ Moreover, the operator $S:\cH(L)\longrightarrow\cH(K)$ determined by
\begin{equation}
\label{TSST}
(SF)(z)=Z(z)F(p(z)),\quad F\in\cH(L), z\in p^{-1}(\Omega_0),
\end{equation}
is unitary. 
\end{Pn}
\begin{proof}
Consider a linear relation in $\cH(K)\times\cH(L)$ spanned by
$$(K(\cdot,w),L(\cdot,p(w))Z(w)^*),\quad w\in p^{-1}(\Omega_0).$$
Since 
$$\|K(\cdot,w)\|^2_{\cH(K)}=K(w,w)=Z(w)L(p(w),p(w))Z(w)^*=\|L(\cdot,p(w))Z(w)^*\|^2_{\cH(L)}$$
and since $\spa\{K_w\,:\,w\in p^{-1}(\Omega_0)\}$ is dense in $\cH(K)$, the above relation is the graph of an isometry
$T\in\bL(\cH(K),\cH(L)).$ The adjoint of $T$ is the operator $S$. In view of Corollary \ref{refmist}, $S$ is injective and hence unitary.

\end{proof}

\begin{Rk}
{\rm In the special case where the kernel $L$ is block diagonal, $L={\rm diag}~(L_1,\ldots, L_N)$, with $L_1,\ldots, L_N$ complex-valued positive definite kernels, we have the orthogonal decomposition
\[
\mathcal H(L)=\oplus_{j=1}^N \mathcal H(L_j),
\]
and, with $S$ as in \eqref{TSST} we can define operators $S_1,\ldots, S_N$ via $S=\begin{pmatrix}S_1&\cdots &S_N\end{pmatrix}$. These operators are given by
\eqref{explicit} and satisfy the Cuntz relations.
 For the theory (and applications) of representations of Cuntz relations by operators in Hilbert space, see e.g. 
\cite{Cun77,MR2277210,MR1743534}.
}
\end{Rk}
{\bf Acknolwedgments:} We thank Professor Vladimir Bolotnikov  for enlightening discussions and for encouragements and  helpful suggestions.
\bibliographystyle{plain}
%\bibliography{/users/faculty/math/dany/Travaux_courants/bib/all}
%\bibliography{/Users/user/Desktop/Travaux_courants/bib/all}
%\bibliography{all}

\begin{thebibliography}{10}

\bibitem{MR2002b:47144}
D.~Alpay.
\newblock {\em The {S}chur algorithm, reproducing kernel spaces and system
  theory}.
\newblock American Mathematical Society, Providence, RI, 2001.
\newblock Translated from the 1998 French original by Stephen S. Wilson,
  Panoramas et Synth\`eses.

\bibitem{CAPB}
D.~Alpay.
\newblock {\em A complex analysis problem book}.
\newblock Birkh\"auser/Springer Basel AG, Basel, 2011.



\bibitem{MR99h:47021}
D.~Alpay, V.~Bolotnikov, and L.~Rodman.
\newblock Tangential interpolation with symmetries and two-point interpolation
  problem for matrix-valued ${H}\sb 2$-functions.
\newblock {\em Integral Equations Operator Theory}, 32(1):1--28, 1998.

\bibitem{MR3050315}
D.~Alpay, P.~Jorgensen, and I.~Lewkowicz.
\newblock Extending wavelet filters: infinite dimensions, the nonrational case,
  and indefinite inner product spaces.
\newblock In {\em Excursions in harmonic analysis. {V}olume 2}, Appl. Numer.
  Harmon. Anal., pages 69--111. Birkh\"auser/Springer, New York, 2013.

\bibitem{ajlm2}
D.~Alpay, P.~Jorgensen, I.~Lewkowicz, and I.~Martziano.
\newblock Infinite products representation for kernels and iteration of
  functions.
\newblock Operator Theory: Advances and Applications, vol. 244 (2015), pp. 67-87.

\bibitem{ajlm1}
D.~Alpay, P.~Jorgensen, I.~Lewkowicz, and I.~Marziano.
\newblock {Representation formulas for Hardy space functions through the Cuntz
  relations and new interpolation problems}.
\newblock In Xiaoping Shen and Ahmed Zayed, editors, {\em {Multiscale signal
  analysis and modeling}}, pages 161--182. Springer, 2013.

\bibitem{av3}
D.~Alpay and V.~Vinnikov.
\newblock Finite dimensional de {B}ranges spaces on {R}iemann surfaces.
\newblock {\em J. Funct. Anal.}, 189(2):283--324, 2002.

\bibitem{bgr}
J.~Ball, I.~Gohberg, and L.~Rodman.
\newblock {\em Interpolation of rational matrix functions}, volume~45 of {\em
  Operator {T}heory: {A}dvances and {A}pplications}.
\newblock Birkh{\" a}user Verlag, Basel, 1990.

\bibitem{bgk1}
H.~Bart, I.~Gohberg, and M.A. Kaashoek.
\newblock {\em Minimal factorization of matrix and operator functions},
  volume~1 of {\em {Operator {T}heory: {A}dvances and {A}pplications}}.
\newblock Birkh{\" a}user Verlag, Basel, 1979.

\bibitem{MR2663312}
H. Bart, I. Gohberg, M. Kaashoek, and A. Ran.
\newblock {\em A state space approach to canonical factorization with
  applications}, volume 200 of {\em Operator Theory: Advances and
  Applications}.
\newblock Birkh\"auser Verlag, Basel; Birkh\"auser Verlag, Basel, 2010.
\newblock Linear Operators and Linear Systems.

\bibitem{BrJo02a}
O.~Bratteli and P.~Jorgensen.
\newblock {\em Wavelets through a looking glass}.
\newblock Applied and Numerical Harmonic Analysis. Birkh\"auser Boston Inc.,
  Boston, MA, 2002.


\bibitem{MR1743534}
O. Bratteli, D. Evans, and P. Jorgensen.
\newblock Compactly supported wavelets and representations of the {C}untz
  relations.
\newblock {\em Appl. Comput. Harmon. Anal.}, 8(2):166--196, 2000.

\bibitem{Cun77}
J. Cuntz.
\newblock Simple {$C\sp*$}-algebras generated by isometries.
\newblock {\em Comm. Math. Phys.}, 57(2):173--185, 1977.


\bibitem{Dev13}
R. Devaney.
\newblock Singular perturbations of complex polynomials.
\newblock {\em Bull. Amer. Math. Soc. (N.S.)}, 50(3):391--429, 2013.

\bibitem{MR2945156}
D. Dutkay and P. Jorgensen.
\newblock Spectral measures and {C}untz algebras.
\newblock {\em Math. Comp.}, 81(280):2275--2301, 2012.

\bibitem{FM13}
J.~C. Ferreira and V.~A. Menegatto.
\newblock Positive definiteness, reproducing kernel {H}ilbert spaces and
  beyond.
\newblock {\em Ann. Funct. Anal.}, 4(1):64--88, 2013.

\bibitem{freitag2}
E.~Freitag and R.~Busam.
\newblock {\em {Complex analysis}}.
\newblock Springer, 2005.


\bibitem{MR2277210}
P.~E.~T.  Jorgensen.
\newblock Certain representations of the {C}untz relations, and a question on
  wavelets decompositions.
\newblock In {\em Operator theory, operator algebras, and applications}, volume
  414 of {\em Contemp. Math.}, pages 165--188. Amer. Math. Soc., Providence,
  RI, 2006.

\bibitem{MR2097020}
P.~E.~T.  Jorgensen.
\newblock Iterated function systems, representations, and {H}ilbert space.
\newblock {\em Internat. J. Math.}, 15(8):813--832, 2004.

\bibitem{MR2254502}
P.~E.~T.  Jorgensen.
\newblock {\em Analysis and probability: wavelets, signals, fractals}, volume
  234 of {\em Graduate Texts in Mathematics}.
\newblock Springer, New York, 2006.


\bibitem{MR1284951}
P.~E.~T. Jorgensen, L.~M. Schmitt, and R.~F. Werner.
\newblock {$q$}-relations and stability of {$C^\ast$}-isomorphism classes.
\newblock In {\em Algebraic methods in operator theory}, pages 261--271.
  Birkh\"auser Boston, Boston, MA, 1994.




\bibitem{MR0255260}
R.~E. Kalman, P.~L. Falb, and M.~A. Arbib.
\newblock {\em Topics in mathematical system theory}.
\newblock McGraw-Hill Book Co., New York, 1969.

\bibitem{MBM+14}
B. Maayah, S. Bushnaq, S. Momani, and O. Abu~Arqub.
\newblock Iterative {M}ultistep {R}eproducing {K}ernel {H}ilbert {S}pace
  {M}ethod for {S}olving {S}trongly {N}onlinear {O}scillators.
\newblock {\em Adv. Math. Phys.}, pages Art. ID 758195, 7, 2014.

\bibitem{meschkowski}
H.~Meschkovski.
\newblock {\em Hilbertsche {R\"aume} mit {K}ernfunktion}.
\newblock Springer--{V}erlag, 1962.

\bibitem{Mou12}
M.~Mouattamid.
\newblock Recursive reproducing kernels {H}ilbert spaces using the theory of
  power kernels.
\newblock {\em Anal. Theory Appl.}, 28(2):111--124, 2012.

\bibitem{NZ13}
M. Nashed and Q. Sun.
\newblock Function spaces for sampling expansions.
\newblock In {\em Multiscale signal analysis and modeling}, pages 81--104.
  Springer, New York, 2013.

\bibitem{saitoh}
S.~Saitoh.
\newblock {\em Theory of reproducing kernels and its applications}, volume 189.
\newblock Longman scientific and technical, 1988.

\bibitem{Wu13}
Q. Wu.
\newblock Regularization networks with indefinite kernels.
\newblock {\em J. Approx. Theory}, 166:1--18, 2013.

\end{thebibliography}
\def\cprime{$'$} \def\cprime{$'$} \def\cprime{$'$}
  \def\lfhook#1{\setbox0=\hbox{#1}{\ooalign{\hidewidth
  \lower1.5ex\hbox{'}\hidewidth\crcr\unhbox0}}} \def\cprime{$'$}
  \def\cfgrv#1{\ifmmode\setbox7\hbox{$\accent"5E#1$}\else
  \setbox7\hbox{\accent"5E#1}\penalty 10000\relax\fi\raise 1\ht7
  \hbox{\lower1.05ex\hbox to 1\wd7{\hss\accent"12\hss}}\penalty 10000
  \hskip-1\wd7\penalty 10000\box7} \def\cprime{$'$} \def\cprime{$'$}
  \def\cprime{$'$} \def\cprime{$'$} \def\cprime{$'$}

\end{document}